\input amstex
\input amsppt.sty
\magnification=\magstep1
\hsize=30truecc
\vsize=22.2truecm
\baselineskip=16truept
\NoBlackBoxes
\TagsOnRight \pageno=1 \nologo
\def\Z{\Bbb Z}
\def\N{\Bbb N}

\def\l{\left}
\def\r{\right}
\def\bg{\bigg}
\def\({\bg(}
\def\[{\bg\lfloor}
\def\){\bg)}
\def\]{\bg\rfloor}
\def\t{\text}
\def\f{\frac}

\def\B{B_{k+1}}

\def\p{\ (\roman{mod}\ p)}

\def\sm{\setminus}

\def\bi{\binom}
\def\eq{\equiv}

\def\ls{\leqslant}
\def\gs{\geqslant}
\def\mo{\roman{mod}}

\def\ve{\varepsilon}
\def\al{\alpha}
\def\da{\delta}

\def\Proof{\noindent{\it Proof}}

\def\Remark{\medskip\noindent{\it  Remark}}

\hbox {Sci. China Math. 53(2010), no.\,9, 2473--2488.}
\bigskip
\topmatter
\title Binomial coefficients, Catalan numbers and Lucas quotients\endtitle
\author Zhi-Wei Sun\endauthor
\date {\it Dedicated to Prof. Yuan Wang on the occasion of his 80th
birthday}\enddate
 \leftheadtext{Zhi-Wei Sun} \rightheadtext{Binomial
coefficients, Catalan numbers, Lucas quotients}
\affil Department of Mathematics, Nanjing University\\
 Nanjing 210093, People's Republic of China
  \\  zwsun\@nju.edu.cn
  \\ {\tt http://math.nju.edu.cn/$\sim$zwsun}
\endaffil
\abstract Let $p$ be an odd prime and let $a,m\in\Z$ with $a>0$ and $p\nmid m$.
In this paper we determine $\sum_{k=0}^{p^a-1}\bi{2k}{k+d}/m^k$ mod $p^2$ for $d=0,1$; for example,
$$\sum_{k=0}^{p^a-1}\f{\bi{2k}k}{m^k}\eq\l(\f{m^2-4m}{p^a}\r)+\l(\f{m^2-4m}{p^{a-1}}\r)u_{p-(\f{m^2-4m}{p})}
\ (\mo\ p^2),$$
where $(-)$ is the Jacobi symbol and
 $\{u_n\}_{n\gs0}$ is the Lucas sequence given by $u_0=0,\ u_1=1$ and $u_{n+1}=(m-2)u_n-u_{n-1}$ ($n=1,2,3,\ldots)$.
As an application, we determine $\sum_{0<k<p^a,\,k\eq r\,(\mo\ p-1)}C_k$ modulo $p^2$
for any integer $r$,
where $C_k$ denotes the Catalan number $\bi{2k}k/(k+1)$. We also pose some related conjectures.
\endabstract
\thanks 2010 {\it Mathematics Subject Classification}.\,Primary 11B65;
Secondary 05A10, 11A07, 11B39.
\newline\indent Supported by the National Natural Science
Foundation (grant 10871087) and the Overseas Cooperation Fund (grant 10928101) of China.
\endthanks
\endtopmatter
\document

\heading{1. Introduction}\endheading

  The well-known Catalan numbers are given by
$$C_k=\f1{k+1}\bi{2k}k=\bi{2k}k-\bi{2k}{k+1}\quad (k\in\N=\{0,1,2,\ldots\}).$$
 They have lots of combinatorial interpretations, see, e.g., [St, pp.\,219--229].

Let $p$ be a prime. In 2006 H. Pan and Z. W. Sun [PS]
obtained some congruences involving Catalan numbers; for example, (1.16) in [PS] yields
$$\sum_{k=1}^{p-1}C_k\eq\f32\l(\l(\f p3\r)-1\r)\ (\mo\ p),$$
where $(-)$ is the Jacobi symbol.
In a recent paper [ST1] Sun and Tauraso
investigated $\sum_{k=0}^{p^a-1}\bi{2k}{k+d}/m^k$ and $\sum_{k=1}^{p-1}\bi{2k}{k+d}/(km^{k-1})$ modulo $p$
via Lucas sequences, where $d$ is an integer among $0,\ldots,p^a$ and  $m$ is an integer not divisible by $p$.
By Sun and R. Tauraso [ST2,
Corollary 1.1], for any $a\in\Z^+=\{1,2,3,\ldots\}$ we have
$$\sum_{k=0}^{p^a-1}\bi{2k}k\eq\l(\f{p^a}3\r)\ (\mo\ p^2)\tag1.1$$
and
$$\sum_{k=0}^{p^a-1}\bi{2k}{k+1}\eq\l(\f{p^a-1}3\r)-p\da_{p,3}\ (\mo\ p^2),\tag1.2$$
where the Kronecker symbol $\da_{m,n}$ takes $1$ or $0$ according as $m=n$ or not.

Let $A\in\Z$ and $B\in\Z\sm\{0\}$. The Lucas sequences $u_n=u_n(A,B)\ (n\in\N)$ and $v_n=v_n(A,B)\ (n\in\N)$
are defined as follows:
$$u_0=0,\ u_1=1,\ \t{and}\ u_{n+1}=Au_n-Bu_{n-1}\ (n=1,2,3,\ldots)$$
and
$$v_0=2,\ v_1=A,\ \t{and}\ v_{n+1}=Av_n-Bv_{n-1}\ (n=1,2,3,\ldots).$$
The characteristic equation $x^2-Ax+B=0$ has two roots
$$\al=\f{A+\sqrt{\Delta}}2\quad\t{and}\quad\beta=\f{A-\sqrt{\Delta}}2,$$
where $\Delta=A^2-4B$. By induction, one can easily get the following well-known formulae:
$$(\al-\beta)u_n=\al^n-\beta^n\quad\ \t{and}\ \ v_n=\al^n+\beta^n.$$
In the case $\al=\beta$ (i.e., $\Delta=0$), clearly $u_n=n(A/2)^{n-1}$ for all $n\in\Z^+$.
If $p$ is an odd prime not dividing $B$, then it is known that $p\mid u_{p-(\f{\Delta}p)}$ (see, e.g., [S06]),
and we call the integer $u_{p-(\f{\Delta}p)}/p$ a {\it Lucas quotient}.
There are many congruences for some special Lucas quotients such as Fibonacci quotients and Pell quotients.
(Cf. [SS] and [S02].)

In this paper we establish the following general theorem
which includes some previous congruences as special cases
and relates binomial coefficients to Lucas quotients.

\proclaim{Theorem 1.1} Let $p$ be an odd prime and let $a\in\Z^+$. Let $m$ be any integer not divisible by $p$
and set $\Delta=m(m-4)$.
Then we have
$$\sum_{k=0}^{p^a-1}\f{\bi{2k}k}{m^k}\eq\l(\f{\Delta}{p^a}\r)
+\l(\f{\Delta}{p^{a-1}}\r)u_{p-(\f{\Delta}{p})}(m-2,1)\ (\mo\ p^2)\tag1.3$$
and
$$\aligned&\sum_{k=0}^{p^a-1}\f{\bi{2k}{k+1}}{m^k}
\eq1-m^{p-1}+\l(\f m2-1\r)\(\l(\f{\Delta}{p^a}\r)-1\)
\\&\qquad+\l(\f m2-1\r)\l(\f{\Delta}{p^{a-1}}\r)u_{p-(\f{\Delta}{p})}(m-2,1)\ (\mo\ p^2).
\endaligned\tag1.4$$
Consequently,
$$\aligned&\sum_{k=1}^{p^a-1}\f{\bi{2k+1}k}{m^k}+m^{p-1}-1
\\\eq&\f m2\(\l(\f{\Delta}{p^a}\r)-1+\l(\f{\Delta}{p^{a-1}}\r)u_{p-(\f{\Delta}p)}(m-2,1)\)\ (\mo\ p^2)
\endaligned\tag1.5$$
and
$$\aligned\sum_{k=1}^{p^a-1}\f{C_k}{m^k}
\eq&m^{p-1}-1-\f{m-4}2\(\l(\f{\Delta}{p^a}\r)-1\)
\\&-\f{m-4}2\l(\f{\Delta}{p^{a-1}}\r)u_{p-(\f{\Delta}{p})}(m-2,1)\ (\mo\ p^2).
\endaligned\tag1.6$$
\endproclaim

Here is a consequence of Theorem 1.1.

\proclaim{Corollary 1.1} Let $p$ be an odd prime
 and let $a\in\Z^+$. Then
$$\sum_{k=0}^{p^a-1}\f{\bi{2k}k}{2^k}\eq(-1)^{(p^a-1)/2}\ (\mo\ p^2)
\ \t{and}\ \sum_{k=1}^{p^a-1}\f{\bi{2k}{k+1}}{2^k}\eq1-2^{p-1}\ (\mo\ p^2).$$
Also,
$$\sum_{k=1}^{p^a-1}\f{\bi{2k}{k+1}}{4^k}\eq p\da_{a,1}-4^{p-1}\ (\mo\ p^2)
\ \t{and}\ \sum_{k=1}^{p^a-1}\f{C_k}{4^k}\eq 2^p-2\ (\mo\ p^2).$$
If $p\not=3$ then
$$\sum_{k=0}^{p^a-1}\f{\bi{2k}k}{3^k}\eq\l(\f{p^a}3\r)\ (\mo\ p^2)
\ \t{and}\ \sum_{k=1}^{p^a-1}\f{C_k}{3^k}\eq 3^{p-1}-1+\f{(\f{p^a}3)-1}2\ (\mo\ p^2).$$
When $p\not=5$ we have
$$\align\sum_{k=0}^{p^a-1}(-1)^k\bi{2k}k\eq&\l(\f{p^a}5\r)\l(1-2F_{p-(\f p5)}\r)\ (\mo\ p^2),
\\ \sum_{k=1}^{p^a-1}(-1)^kC_k\eq&\f52\l(\l(\f{p^a}5\r)-1\r)-5\l(\f{p^a}5\r)F_{p-(\f p5)}\ (\mo\ p^2),
\\\sum_{k=0}^{p^a-1}\f{\bi{2k}k}{5^k}\eq&\l(\f{p^a}5\r)\l(1+2F_{p-(\f{p}5)}\r)\ (\mo\ p^2),
\\\sum_{k=1}^{p^a-1}\f{C_k}{5^k}\eq&\f{1-(\f{p^a}5)}2-\l(\f{p^a}5\r)F_{p-(\f p5)}\ (\mo\ p^2),
\endalign$$
where $\{F_n\}_{n\gs0}$ is the well-known Fibonacci sequence defined by
$$F_0=0,\ F_1=1,\ \t{and}\ F_{n+1}=F_n+F_{n-1}\ (n=1,2,3,\ldots).$$
\endproclaim

\Remark\ 1.1. (i) There is a closed formula for the sum $\sum_{k=0}^n\bi{2k}k/4^k$.
In fact, $\bi{2k}k=(-4)^k\bi{-1/2}k$ for $k\in\N$ and hence
$$\sum_{k=0}^n\f{\bi{2k}k}{4^k}=(-1)^n\sum_{k=0}^n\bi{-1}{n-k}\bi{-1/2}k=(-1)^n\bi{-3/2}{n}=\f{2n+1}{4^n}\bi{2n}n$$
by the Chu-Vandermonde identity
$$\sum_{k=0}^n\bi xk\bi y{n-k}=\bi{x+y}n$$
(see, e.g., [GKP, p.\,169]).

(ii) In [ST1] the authors conjectured that if $p\not=2,5$ is a prime
and $a\in\Z^+$ then
$$\sum_{k=0}^{p^a-1}(-1)^k\bi{2k}k\eq\l(\f{p^a}5\r)\l(1-2F_{p^a-(\f {p^a}5)}\r)\ (\mo\ p^3).$$
(Note that $F_{p^a-(\f{p^a}5)}\eq F_{p-(\f p5)}\ (\mo\ p^2)$ by Lemma 2.3.)
This seems difficult.
Those primes $p>5$ satisfying $p^2\mid F_{p-(\f p5)}$ are called Wall-Sun-Sun primes (cf. [CP, p.\,32]).
Up to now none of this kind of primes has been found though
it is conjectured that there should be infinitely many Wall-Sun-Sun primes.
\medskip

By Corollary 1.1, if $p$ is an odd prime then
$$\sum_{k=0}^{p-1}\f{\bi{2k}k}{2^k}\eq(-1)^{(p-1)/2}\ (\mo\ p^2).$$
This seems to be a new characterization of odd primes and we have verified
our following conjecture for $n<10^4$ via {\tt Mathematica}.

\proclaim{Conjecture 1.1} If an odd integer $n>1$ satisfies the congruence
$$\sum_{k=0}^{n-1}\f{\bi{2k}k}{2^k}\eq(-1)^{(n-1)/2}\ (\mo\ n^2),$$
then $n$ must be a prime.
\endproclaim

As an application of Theorem 1.1, we will determine the sums
$$\sum\Sb 0<k<p^a\\k\eq r\,(\mo\ p-1)\endSb \bi{2k}{k},\
\sum\Sb 0<k<p^a\\k\eq r\,(\mo\ p-1)\endSb \bi{2k}{k+1},\ \sum\Sb 0<k<p^a\\k\eq r\,(\mo\ p-1)\endSb C_k
$$
modulo $p^2$ for any prime $p$ and integers $a>0$ and $r$.
By (1.1) and (1.2), for $d=0,1$ we have
$$\sum_{k=0}^{p^a-1}\bi{2k}{k+d}\eq\l(\f{p^a-d}3\r)-p\da_{d,1}\da_{p,3}\ (\mo\ p^2).$$
Thus the task for $p=2$ is easy; for example,
$$\align\sum\Sb 0<k<2^a\\k\eq r\,(\mo\ 2-1)\endSb C_k=
&\sum_{k=1}^{2^a-1}C_k=\sum_{k=1}^{2^a-1}\bi{2k}k-\sum_{k=1}^{2^a-1}\bi{2k}{k+1}
\\\eq&\l(\f{2^a}3\r)-1-\l(\f{2^a-1}3\r)\eq\cases1\ (\mo\ 2^2)&\t{if}\ 2\nmid a,\\0\ (\mo\ 2^2)&\t{if}\ 2\mid a.
\endcases
\endalign$$
So we will only handle the main case $p\not=2$.

\proclaim{Theorem 1.2} Let $p$ be an odd prime and let $a\in\Z^+$.

{\rm (i)} If $a$ is odd and $r\in\{1,\ldots,p-1\}$, then
$$\sum\Sb 0<k<p^a\\k\eq r\,(\mo\ p-1)\endSb\bi{2k}{k+d}\eq\bi{2r}{r+d}\ (\mo\ p^2)\quad\t{for}\ d=0,1,\tag1.7$$
and also
$$\sum\Sb 0<k<p^a\\k\eq r\,(\mo\ p-1)\endSb C_k\eq C_r\ (\mo\
p^2).\tag1.8$$

{\rm (ii)} Suppose that $a$ is even. Then, for $r=1,\ldots,p$ we have
$$\sum\Sb 0<k<p^a\\k\eq r\,(\mo\ p-1)\endSb\bi{2k}k\eq
4^r\l(1+\f p2+r(2^{p-1}-1)\r)-pR_p(r)\ (\mo\ p^2),\tag1.9$$
where
$$R_p(r)=\cases\sum_{s=0}^{(p-1)/2-r}\bi{2r+2s}{r+s}/((2s+1)\bi{2s}s)\ &\t{if}\ 0<r\ls (p-1)/2,
\\0&\t{otherwise}.\endcases$$
Also, if $r\in\{1,\ldots,p-1\}$ then
$$\aligned \sum\Sb 0<k<p^a\\k\eq r\,(\mo\ p-1)\endSb\bi{2k}{k+1}
\eq&4^r\l(1+\f p2+(r+2)(2^{p-1}-1)\r)
\\&+p\l(R_p(r)-\f{R_p(r+1)}2\r)\ (\mo\ p^2)
\endaligned\tag1.10$$
and
$$\sum\Sb 0<k<p^a\\k\eq r\,(\mo\ p-1)\endSb C_k
\eq4^r(2-2^p)-p\l(2R_p(r)-\f{R_p(r+1)}2\r)\ (\mo\ p^2).\tag1.11$$
In particular,
$$\sum\Sb 0<k<p^a\\k\eq r\,(\mo\ p-1)\endSb C_k
\eq4^r(2-2^p)\ (\mo\ p^2)\ \ \t{for}\ r=\f{p+1}2,\ldots,p-1.\tag1.12$$
\endproclaim
\Remark\ 1.2. If $p$ is an odd prime and $a\in\Z^+$ is even, then by (1.11) we have
$$\sum\Sb 0<k<p^a\\k\eq r\,(\mo\ p-1)\endSb C_k\eq0\ (\mo\ p)\quad\t{for all}\ r\in\Z.$$
The author would like to see any combinatorial interpretation for this.
\medskip

\proclaim{Corollary 1.2} Let $p$ be an odd prime and let $a\in\Z^+$. Then
$$\sum\Sb 0<k<p^a\\k\eq0\,(\mo\ p-1)\endSb C_k\eq \cases-2p-1\ (\mo\ p^2)&\t{if}\ 2\nmid a,
\\2-2^p\ (\mo\ p^2)&\t{if}\ 2\mid a;\endcases\tag1.13$$
$$\sum\Sb 0<k<p^a\\k\eq1\ (\mo\ p-1)\endSb C_k\eq\cases1\ (\mo\ p^2)
&\t{if}\ 2\nmid a,\\4(2-2^p)+2p\ (\mo\ p^2)&\t{if}\ 2\mid a;
\endcases\tag1.14$$
and
$$\sum\Sb 0<k<p^a\\k\eq(p-1)/2\,(\mo\ p-1)\endSb C_k
\eq\cases(-1)^{(p-1)/2}2(2^p-p-1)\ (\mo\ p^2)&\t{if}\ 2\nmid a,
\\2-2^p+(-1)^{(p+1)/2}2p\ (\mo\ p^2)&\t{if}\ 2\mid a.
\endcases\tag1.15$$
\endproclaim

Now we pose some new conjectures.

\proclaim{Conjecture 1.2} Let $p$ be any prime and let $r$ be an
integer. For $a\in\N$ define
$$S_r(p^a)=\sum\Sb 0<k<p^a\\k\eq r\,(\mo\ p-1)\endSb C_k.$$
Then, for any $a\in\N$ we have
$$S_r(p^{a+2})\eq S_r(p^a)\ (\mo\ p^{(1+\da_{p,2})(a+1)}).$$
Furthermore,
$$\f{S_r(p^{a+2})-S_r(p^a)}{p^{(1+\da_{p,2})(a+1)}}+p(\da_{p^a,2}+\da_{p^a,3}) \ \ \mo\ p^2$$
does not depend on $a\in\Z^+$.
\endproclaim

\proclaim{Conjecture 1.3} Let $p$ be a prime, and let
$d\in\{0,\ldots,p\}$ and $r\in\Z$. For $a\in\N$ define
$$T^{(d)}_r(p^a)=\sum\Sb
0<k<p^a\\k\eq r\,(\mo\ p-1)\endSb\bi{2k}{k+d}.$$
Then, for any $a\in\N$ we have
$$T^{(d)}_r(p^{a+2})\eq T^{(d)}_r(p^a)\ (\mo\ p^a);$$
furthermore
$$\f{T^{(d)}_r(p^{a+2})-T^{(d)}_r(p^a)}{p^a}\ \mo\ p$$
does not depend on $a\in\Z^+$.
If $a\in\N$ and $d<p=2$, then
$$T^{(d)}_r(2^{a+2})\eq T^{(d)}_r(2^a)\ (\mo\ 2^{2a+2+\delta_{d,0}(1-\da_{a,0})}).$$
If $a\in\Z^+$, $d\in\{0,1\}$ and $p=3$, then
$$T^{(d)}_r(3^{a+2})\eq T^{(d)}_r(3^a)\ (\mo\ 3^{a+1+\delta_{d,1}(1-\da_{a,1})}).$$
\endproclaim

Given a positive integer $h$, two kinds of Catalan numbers of order $h$ are defined as follows:
$$C_k^{(h)}=\f1{hk+1}\bi{(h+1)k}k=\bi{(h+1)k}k-h\bi{(h+1)k}{k-1}\ \ (k\in\N)$$
and
$$\bar C_k^{(h)}=\f h{k+1}\bi{(h+1)k}k=h\bi{(h+1)k}k-\bi{(h+1)k}{k+1}\ \ (k\in\N).$$
In [ZPS] and [S09], the authors gave various congruences involving higher-order Catalan numbers.
In particular, Sun [S09] proved that for any prime $p>3$ and $a\in\Z^+$ with $6\mid a$ we have
the congruence
 $$\sum\Sb 0<k<p^a\\k\eq r\,(\mo\ p-1)\endSb\bi{3k}{k+d}\eq2^{d+3-2r}3^{3r-2}\ (\mo\ p)$$
for all $d\in\{0,\pm1\}$ and $r\in\Z$; consequently,
$$\sum\Sb 0<k<p^a\\k\eq r\,(\mo\ p-1)\endSb C_k^{(2)}\eq \sum\Sb 0<k<p^a\\k\eq r\,
(\mo\ p-1)\endSb \bar C_k^{(2)}\eq0\ (\mo\ p)$$
for any $r\in\Z$.

Here is our conjecture involving Catalan numbers of order 2.

\proclaim{Conjecture 1.4} Let $p$ be any prime, and set
  $$C(p^a)=\sum\Sb 0<k<p^a\\k\eq0\,(\mo\ p-1)\endSb C_k^{(2)}\ \t{and}
  \ \bar C(p^a)=\sum\Sb 0<k<p^a\\k\eq0\,(\mo\ p-1)\endSb \bar C_k^{(2)}\quad\t{for}\ a\in\Z^+.$$
Then we have
$$C(p^a)\eq\cases0\ (\mo\ p)&\t{if}\ a\eq0\ (\mo\ 6),
\\\da_{p,2}\ (\mo\ p)&\t{if}\ a\eq1\ (\mo\ 6),
\\-((\f p3)+1)/2\ (\mo\ p)&\t{if}\ a\eq2\ (\mo\ 6),
\\((\f p3)-1)/2+\da_{p,2}\ (\mo\ p)&\t{if}\ a\eq 3\ (\mo\ 6),
\\(1-(\f p3))/2\ (\mo\ p)&\t{if}\ a\eq4\ (\mo\ 6),
\\\da_{p,2}-1\ (\mo\ p)&\t{if}\ a\eq5\ (\mo\ 6);\endcases$$
and
$$\bar C(p^a)\eq\cases0\ (\mo\ p)&\t{if}\ a\eq0\ (\mo\ 6),
\\-2+\da_{p,2}\ (\mo\ p)&\t{if}\ a\eq\pm1\ (\mo\ 6),
\\-1-2(\f p3)\ (\mo\ p)&\t{if}\ a\eq\pm2\ (\mo\ 6),
\\2(\f p3)-1+\da_{p,2}\ (\mo\ p)&\t{if}\ a\eq 3\ (\mo\ 6).
\endcases$$
\endproclaim

We will prove Theorem 1.1 and Corollary 1.1 in Section 2,
 and  show Theorem 1.2 and Corollary 1.2 in Section 3.

\heading{2. Proof of Theorem 1.1}\endheading

\proclaim{Lemma 2.1} Let $p$ be a prime and let $a,m\in\Z$ with $a>0$ and $p\nmid m$. Then
$$\sum_{k=1}^{p^a-1}\f{\bi{2k}{k+1}}{m^k}+(m^{p-1}-1)\eq
\f{m-2}2\sum_{k=1}^{p^a-1}\f{\bi{2k}k}{m^k}+p\da_{p,2}\ (\mo\ p^2).
\tag2.1$$
\endproclaim
\Proof.
Observe that
$$\align\sum_{k=0}^{p^a-1}\f{\bi{2k}k+\bi{2k}{k+1}}{m^k}
=&\f12\sum_{k=0}^{p^a-1}\f{\bi{2(k+1)}{k+1}}{m^k}=\f12\sum_{k=1}^{p^a}\f{\bi{2k}k}{m^{k-1}}
\\=&\f12\(\sum_{k=0}^{p^a-1}\f{\bi{2k}k}{m^{k-1}}-m+\f{\bi{2p^a}{p^a}}{m^{p^a-1}}\)
\\=&\f m2\sum_{k=0}^{p^a-1}\f{\bi{2k}k}{m^k}-\f m2+\f{\bi{2p^a-1}{p^a-1}}{m^{p^a-1}}.
\endalign$$
Clearly we have
$$\align&\bi{2p^a-1}{p^a-1}=\prod_{k=1}^{p^a-1}\l(1+\f{p^a}k\r)
\\\eq&1+\f12\sum_{k=1}^{p^a-1}\(\f {p^a}k+\f{p^a}{p^a-k}\)
\eq1+p\da_{p,2}\ (\mo\ p^2).\endalign$$
(See also [ST2, Lemma 2.2].)
Note that
$$\f1{m^{p^a-1}}\eq\f1{m^{p-1}}\eq2-m^{p-1}\ (\mo\ p^2)$$
since $m^{p(p-1)}\eq1\ (\mo\ p^2)$ and $(m^{p-1}-1)^2\eq0\ (\mo\ p^2)$ by Euler's theorem and Fermat's little theorem.
Therefore
$$\align\sum_{k=1}^{p^a-1}\f{\bi{2k}{k+1}}{m^k}\eq&\l(\f m2-1\r)\sum_{k=1}^{p^a-1}\f{\bi{2k}k}{m^k}
+1-m^{p-1}+p\da_{p,2}(2-m^{p-1})
\\\eq&\f{m-2}2\sum_{k=1}^{p^a-1}\f{\bi{2k}k}{m^k}+1-m^{p-1}+p\da_{p,2}\ (\mo\ p^2).
\endalign$$
This concludes the proof. \qed

\proclaim{Lemma 2.2} Let $p$ be any prime and let $a\in\Z^+$. Let $m$ ba an integer not divisible by $p$.
Then
$$\f{m^{p-1}}2\sum_{k=0}^{p^a-1}\f{\bi{2k}k}{m^k}+\f{u_{p^a}(m-2,1)}2\eq u_{p^a}(m,m)\ (\mo\ p^2).\tag2.2$$
\endproclaim
\Proof. By [ST1, Theorem 2.1],
$$\sum_{k=0}^{p^a-1}\bi{2k}km^{p^a-1-k}=\sum_{k=0}^{p^a-1}\bi{2p^a}ku_{p^a-k}(m-2,1).$$
For $k\in\{1,\ldots,p^a-1\}$, clearly
$$\bi{2p^a-1}{k-1}=\prod_{0<j<k}\f{2p^a-j}j\eq\prod_{0<j<k}\f{p^a-j}j=\bi{p^a-1}{k-1}\ (\mo\ p)$$
and hence
$$\f12\bi{2p^a}k=\f{p^a}k\bi{2p^a-1}{k-1}\eq\f{p^a}k\bi{p^a-1}{k-1}=\bi{p^a}k\ (\mo\ p^2).$$
Therefore
$$\align&\f{m^{p^a-1}}2\sum_{k=0}^{p^a-1}\f{\bi{2k}k}{m^k}+\f{u_{p^a}(m-2,1)}2
\\=&\f12\sum_{k=1}^{p^a-1}\bi{2p^a}ku_{p^a-k}(m-2,1)+u_{p^a}(m-2,1)
\\\eq&\sum_{k=1}^{p^a}\bi{p^a}ku_{p^a-k}(m-2,1)+u_{p^a}(m-2,1)
\\\eq&\sum_{j=0}^{p^a}\bi{p^a}ju_j(m-2,1)\ (\mo\ p^2).
\endalign$$

If $\Delta=(m-2)^2-4=m^2-4m\not=0$ then
$$\align&\sum_{j=0}^{p^a}\bi{p^a}ju_j(m-2,1)
\\=&\sum_{j=0}^{p^a}\bi{p^a}j\f1{\sqrt{\Delta}}\l(\l(\f{m-2+\sqrt{\Delta}}2\r)^j-\l(\f{m-2-\sqrt{\Delta}}2\r)^j\r)
\\=&\f1{\sqrt{\Delta}}\l(\l(\f{m+\sqrt{\Delta}}2\r)^{p^a}-\l(\f{m-\sqrt{\Delta}}2\r)^{p^a}\r)=u_{p^a}(m,m).
\endalign$$
In the case $\Delta=0$ (i.e., $m=4$), we have
$$\sum_{j=0}^{p^a}\bi{p^a}ju_j(2,1)=\sum_{j=0}^{p^a}
\bi{p^a}jj=p^a\sum_{j=1}^{p^a}\bi{p^a-1}{j-1}=p^a2^{p^a-1}=u_{p^a}(4,4).$$

In view of the above, it suffices to show that
$$\f{m^{p^a-1}-m^{p-1}}2\eq0\ (\mo\ p^2).\tag2.3$$
This follows from Euler's theorem when $p\not=2$.
If $p=2$, then (2.3) holds since $2\nmid m$ and $m^p=m^2\eq1\ (\mo\ 2^3)$.
We are done. \qed

\medskip

Now we need a lemma on Lucas sequences.

\proclaim{Lemma 2.3} Let $p$ be a prime, and let $a\in\Z^+$ and $A,B\in\Z$. Then
$$v_{p^a}(A,B)\eq v_{p^{a-1}}(A,B)\ (\mo\ p^a).\tag2.4$$
If $p\not=2$, then
$$u_{p^a}(A,B)\eq\l(\f{\Delta}p\r)u_{p^{a-1}}(A,B)\ (\mo\ p^a),\tag2.5$$
where $\Delta=A^2-4B$.
When $p\nmid 2B\Delta$, we have
$$u_{p^a-(\f{\Delta}{p^a})}(A,B)\eq\cases
B^{((\f{\Delta}{p^{a-1}})-(\f{\Delta}{p^a}))/2}(\f{\Delta}p)u_{p^{a-1}-(\f{\Delta}{p^{a-1}})}(A,B)
\ (\mo\ p^a)\\
B^{((\f{\Delta}p)-(\f{\Delta}{p^a}))/2}
(\f{\Delta}{p^{a-1}})u_{p-(\f{\Delta}p)}(A,B)
\ (\mo\ p^2)
\\0\ (\mo\ p)
\endcases\tag2.6$$
and
$$v_{p^a-(\f{\Delta}{p^a})}(A,B)\eq\cases
B^{((\f{\Delta}{p^{a-1}})-(\f{\Delta}{p^a}))/2}v_{p^{a-1}-(\f{\Delta}{p^{a-1}})}(A,B)\ (\mo\ p^a)\\
\\B^{((\f{\Delta}p)-(\f{\Delta}{p^a}))/2}v_{p-(\f{\Delta}p)}(A,B)\ (\mo\ p^2)
\\2B^{(1-(\f{\Delta}{p^a}))/2}\ (\mo\ p).
\endcases\tag2.7$$
\endproclaim
\Proof. For convenience we let $u_n=u_n(A,B)$ and $v_n=v_n(A,B)$ for all $n\in\N$.
We split our proof into several steps.

(i) By a known result of W. J\"anichen [J] (see also [Sm] and [V]),
if $\prod_{j=1}^m(x-\al_j)\in\Z[x]$ then
$$\al_1^{p^a}+\cdots+\al_m^{p^a}\eq \al_1^{p^{a-1}}+\cdots+\al_m^{p^{a-1}}\ (\mo\ p^a).$$
Thus
$$v_{p^a}=\al^{p^a}+\beta^{p^a}\eq \al^{p^{a-1}}+\beta^{p^{a-1}}=v_{p^{a-1}}\ (\mo\ p^a),$$
where $\al$ and $\beta$ be the two roots of the equation $x^2-Ax+B=0$ in the complex field.

(ii) Now we prove that  $p^a\mid u_{p^a}$ under the condition $p\mid \Delta$.

If $\Delta=0$ (i.e., $\al=\beta$), then $A$ is even and $u_n=n(A/2)^{n-1}$ for all $n\in\Z^+$,
in particular $u_{p^a}\eq0\ (\mo\ p^a)$.

Assume $\Delta\not=0$. If $p\not=2$, then
$$u_p\eq u_p\l(A,\f{A^2}4\r)=p\l(\f A2\r)^{p-1}\eq0\ (\mo\ p).$$
When $p=2$, we have $2\mid A$ since $p\mid\Delta$, hence
$u_2=A\eq0\ (\mo\ 2)$. So we always have $p\mid u_p$.
Observe that
$$\align u_{p^{a+1}}=&\f{\al^{p^{a+1}}-\beta^{p^{a+1}}}{\al-\beta}
\\=&\f{\al^{p^a}-\beta^{p^a}}{\al-\beta}\sum_{k=0}^{p-1}(\al^{p^a})^k(\beta^{p^a})^{p-1-k}
\\=&u_{p^a}\sum_{k=0}^{p-1}(\al^k\beta^{p-1-k})^{p^a}
\endalign$$
and
$$\align\sum_{k=0}^{p-1}(\al^k\beta^{p-1-k})^{p^a}\eq&\(\sum_{k=0}^{p-1}\al^k\beta^{p-1-k}\)^{p^a}
\\\eq&\l(\f{\al^p-\beta^p}{\al-\beta}\r)^{p^a}=u_p^{p^a}\eq0\ (\mo\ p).
\endalign$$
Thus, if $p^a\mid u_{p^a}$ then $p^{a+1}\mid u_{p^{a+1}}$.
This concludes our induction proof of the desired congruence $u_{p^a}\eq0\ (\mo\ p^a)$.

(iii) Suppose $p\not=2$. Now we show that
$$u_{p^a}\eq\l(\f{\Delta}{p^a}\r)\ (\mo\ p).$$
By part (ii), this holds when $p\mid \Delta$. In the case $p\nmid\Delta$, since
$$\Delta u_{p^a}=(\al-\beta)^2u_{p^a}=(\al-\beta)(\al^{p^a}-\beta^{p^a})\eq(\al-\beta)^{p^a+1}=\Delta^{(p^a+1)/2}\ (\mo\ p),$$
we have
$$u_{p^a}\eq \Delta^{(p^a-1)/2}\eq\l(\f{\Delta}p\r)^{\sum_{i=0}^{a-1}p^i}
=\l(\f{\Delta}p\r)^a=\l(\f{\Delta}{p^a}\r)\ (\mo\ p).$$

(iv) Assume that $p\not=2$. By part (ii), (2.5) holds when $p\mid\Delta$.
Suppose $p\nmid \Delta$. In view of part (iii),
$$u_{p^a}+\l(\f{\Delta}p\r)u_{p^{a-1}}\eq  2\l(\f{\Delta}{p^a}\r)\not\eq0\ (\mo\ p).$$
For any $n\in\N$ we have
$$v_n^2-\Delta u_n^2=(\al^n+\beta^n)^2-(\al^n-\beta^n)^2=4(\al\beta)^n=4B^n.$$
Thus
$$\Delta (u_{p^a}^2-u_{p^{a-1}}^2)=v_{p^a}^2-4B^{p^a}-(v_{p^{a-1}}^2-4B^{p^{a-1}})$$
and hence
$$\align&\Delta \l(u_{p^a}+\l(\f{\Delta}p\r)u_{p^{a-1}}\r) \l(u_{p^a}-\l(\f{\Delta}p\r)u_{p^{a-1}}\r)
\\=&\l(v_{p^a}+v_{p^{a-1}}\r)\l(v_{p^a}-v_{p^{a-1}}\r)-4(B^{p^a}-B^{p^{a-1}})
\\\eq&0\ (\mo\ p^a)\ \ (\t{by (2.4) and Euler's theorem}).
\endalign$$
So (2.5) follows, for,
$\Delta(u_{p^a}+(\f{\Delta}p)u_{p^{a-1}})$
is relatively prime to $p$.

(v) By induction, for $\ve\in\{\pm1\}$ and $n\in\Z^+$ we have
$$Au_n+\ve v_n=2B^{(1-\ve)/2}u_{n+\ve}\ \ \t{and}\ \ Av_n+\ve \Delta u_n=2B^{(1-\ve)/2}v_{n+\ve}.\tag2.8$$
Therefore, if $p\nmid 2B\Delta$ then
$$\align &u_{p^a-(\f{\Delta}{p^a})}=\f{Au_{p^a}-(\f{\Delta}{p^a})v_{p^a}}{2B^{(1+(\f{\Delta}{p^a}))/2}}
\\\eq&\f{A(\f{\Delta}p)u_{p^{a-1}}-(\f{\Delta}{p^a})v_{p^{a-1}}}{2B^{(1+(\f{\Delta}{p^a}))/2}}
\\\eq&\l(\f{\Delta}p\r)B^{((\f{\Delta}{p^{a-1}}-(\f{\Delta}{p^a}))/2}u_{p^{a-1}-(\f{\Delta}{p^{a-1}})}\ (\mo\ p^a)
\endalign$$
and
$$\align &v_{p^a-(\f{\Delta}{p^a})}=\f{Av_{p^a}-(\f{\Delta}{p^a})\Delta u_{p^a}}{2B^{(1+(\f{\Delta}{p^a}))/2}}
\\\eq&\f{Av_{p^{a-1}}-(\f{\Delta}{p^{a-1}})\Delta u_{p^{a-1}}}{2B^{(1+(\f{\Delta}{p^a}))/2}}
=B^{((\f{\Delta}{p^{a-1}}-(\f{\Delta}{p^a}))/2}v_{p^{a-1}-(\f{\Delta}{p^{a-1}})}\ (\mo\ p^a).
\endalign$$
Note that $u_{p^0-(\f{\Delta}{p^0})}=u_0=0$ and $v_{p^0-(\f{\Delta}{p^0})}=v_0=2$. So both (2.6) and (2.7) hold when
$p\nmid 2B\Delta$.

So far we have completed the proof of Lemma 2.3. \qed

Using Lemma 2.3 we can deduce the following result.
\proclaim{Lemma 2.4} Let $p$ be an odd prime, and let $a,m\in\Z$ with $a>0$ and $p\nmid m$.
Set $\Delta=m^2-4m$. Then
$$2u_{p^a}(m,m)-u_{p^a}(m-2,1)\eq \l(\f{\Delta}{p^a}\r)m^{p-1}+u_{p^a-(\f{\Delta}{p^a})}(m-2,1)\ (\mo\ p^2).\tag2.9$$
\endproclaim
\Proof. By Lemma 2,3,
$$2u_{p^a}(m,m)-u_{p^a}(m-2,1)\eq\l(\f{\Delta}{p^{a-1}}\r)(2u_p(m,m)-u_p(m-2,1))\ (\mo\ p^2)$$
and
$$u_{p^a-(\f{\Delta}{p^a})}(m-2,1)\eq\l(\f{\Delta}{p^{a-1}}\r)u_{p-(\f{\Delta}p)}(m-2,1)\ (\mo\ p^2).$$
So, it suffices to prove (2.9) in the case $a=1$.

Let $\al$ and $\beta$ be the two roots of the equation $x^2-mx+m=0$.
Clearly $(\al-1)+(\beta-1)=m-2$ and $(\al-1)(\beta-1)=1$.
Recall that $\Delta=m^2-4m=(m-2)^2-4$.
If $\Delta\not=0$, then $\al\not=\beta$ and hence
$$\align u_n(m-2,1)=&\f{(\al-1)^n-(\beta-1)^n}{(\al-1)-(\beta-1)}=\f{(\al^2/m)^n-(\beta^2/m)^n}{\al-\beta}
\\=&\f{\al^n-\beta^n}{\al-\beta}\cdot\f{\al^n+\beta^n}{m^n}=\f{u_n(m,m)v_n(m,m)}{m^n}
\endalign$$
for all $n\in\N$. In the case $\Delta=0$ (i.e., $m=4$), as
$$u_n(2,1)=n,\ u_n(4,4)=n2^{n-1}\ \t{and}\ v_n(4,4)=2^{n+1},$$ we also have
$$u_n(m-2,1)=n=\f{n2^{n-1}2^{n+1}}{4^n}=\f{u_n(m,m)v_n(m,m)}{m^n}.$$
So, for any $n\in\N$ we always have
$$u_n(m-2,1)=\f{u_n(m,m)v_n(m,m)}{m^n}.\tag2.10$$

Note that  $v_p(m,m)\eq v_{p^0}(m,m)=m\ (\mo\ p)$ by (2.4).
In view of (2.10) and Lemma 2.3,
$$\align&2u_{p}(m,m)-u_{p}(m-2,1)
\\=&\f{u_{p}(m,m)}{m^{p}}\l(m^{p}-v_{p}(m,m)\r)+u_{p}(m,m)
\\\eq&\f{(\f{\Delta}{p})}m\l(m^{p}-v_{p}(m,m)\r)+u_{p}(m,m)
\\\eq&\l(\f{\Delta}{p}\r)m^{p-1}+u_{p}(m,m)-\l(\f{\Delta}{p}\r)\f{v_{p}(m,m)}m\ (\mo\ p^2).
\endalign$$
Thus, by the above, it suffices to prove the congruence
$$u_{p-(\f{\Delta}{p})}(m,m)\f{v_{p-(\f{\Delta}{p})}(m,m)}{m^{p-(\f{\Delta}{p})}}
\eq u_{p}(m,m)-\l(\f{\Delta}{p}\r)\f{v_{p}(m,m)}m\ (\mo\ p^2).\tag2.11$$

Clearly, $u_{p-(\f{\Delta}p)}(m,m)\eq0\ (\mo\ p)$ by Lemma 2.3. If $p\mid \Delta$ then
$$v_{p-(\f{\Delta}p)}(m,m)=v_p(m,m)\eq m\eq m^p=m^{p-(\f{\Delta}p)}\ (\mo\ p)$$
and hence (2.11) holds.

Now assume that $p\nmid\Delta$. Obviously,
$$\f{v_{p-(\f{\Delta}p)}(m,m)}{m^{p-(\f{\Delta}p)}}\eq \f{2m^{(1-(\f{\Delta}p))/2}}{m^{1-(\f{\Delta}p)}}
=2m^{((\f{\Delta}p)-1)/2}\ (\mo\ p)$$
by (2.7), and
$$u_{p-(\f{\Delta}p)}(m,m)=\f{mu_p(m,m)-(\f{\Delta}p)v_p(m,m)}{2m^{(1+(\f{\Delta}p))/2}}$$
by (2.8). Therefore the left-hand side of (2.11) is congruent to
$$\f{mu_p(m,m)-(\f{\Delta}p)v_p(m,m)}m=u_p(m,m)-\l(\f{\Delta}p\r)\f{v_p(m,m)}m$$
modulo $p^2$. So (2.11) is valid and we are done.
\qed

\medskip
\noindent{\it Proof of Theorem 1.1}.  Clearly (1.3) plus or minus (1.4) yields (1.5) or (1.6).
Also, (1.4) follows from (1.3) by Lemma 2.1.
So, it suffices to prove (1.3).

 Combining Lemmas 2.2--2.4, we get
 $$\align m^{p-1}\sum_{k=0}^{p^a-1}\f{\bi{2k}k}{m^k}
 \eq&2u_{p^a}(m,m)-u_{p^a}(m-2,1)
 \\\eq&\l(\f{\Delta}{p^{a}}\r)m^{p-1}+u_{p^a-(\f{\Delta}{p^a})}(m-2,1)
 \\\eq&\l(\f{\Delta}{p^{a-1}}\r)m^{p-1}\(\l(\f{\Delta}{p}\r)+u_{p-(\f{\Delta}{p})}(m-2,1)\)\ (\mo\ p^2).
 \endalign$$
Therefore (1.3) holds. This concludes the proof. \qed

\medskip
\noindent{\it Proof of Corollary 1.1}. By induction, $u_{2n}(0,1)=0$ and $u_n(2,1)=n$ for all $n\in\N$.
Note also that
$$(-1)^{n-1}u_n(1,1)=u_n(-1,1)=\l(\f n3\r)$$
and
$$(-1)^{n-1}u_n(-3,1)=u_n(3,1)=F_{2n}=F_nL_n,$$
where $L_n=v_n(1,-1)$. By [SS, Corollary 1] (or the proof of Corollary 1.3 of [ST1]), if $p\not=2,5$ then
$L_{p-(\f p5)}\eq 2\l(\f p5\r)\ (\mo\ p^2).$

In view of the above, we can easily deduce the  congruences in Corollary 1.1 by applying Theorem 1.1. \qed

 \heading{3. Proof of Theorem 1.2}\endheading

\proclaim{Lemma 3.1} Let $p$ be an odd prime and let $k\in\Z$. Then
$$\sum_{m=1}^{p-1}{m^{pk}}\eq\cases p-1\ (\mo\ p^2)&\t{if}\ p-1\mid k,\\0\ (\mo\ p^2)&\t{otherwise}.
\endcases$$
\endproclaim
\Proof. For $b,c\in\Z$ clearly $(b+cp)^p\eq b^p\ (\mo\ p^2)$. If $p-1\mid k$, then $m^{pk}\eq1\ (\mo\ p^2)$
by Euler's theorem, and hence $\sum_{m=1}^{p-1}m^{pk}\eq(p-1)\ (\mo\ p^2)$.

Now suppose that $p-1\nmid k$ and let $g$ be a primitive root modulo $p$. Then
$$g^{pk}\sum_{m=1}^{p-1}{m^{pk}}=\sum_{m=1}^{p-1}{(gm)^{pk}}\eq\sum_{r=1}^{p-1}{r^{pk}}\ (\mo\ p^2)$$
and hence
$$(g^{pk}-1)\sum_{m=1}^{p-1}{m^{pk}}\eq0\ (\mo\ p^2).$$
Since $g^{pk}-1$ is not divisible by $p$, we must have
$$\sum_{m=1}^{p-1}{m^{pk}}\eq0\ (\mo\ p^2).$$
This concludes the proof. \qed

\proclaim{Lemma 3.2} Let $p$ be an odd prime and let $a\in\Z^+$. Then, for any $r\in\Z$, we have
$$\sum\Sb 0<k<p^a\\k\eq r\,(\mo\ p-1)\endSb\bi{2k}{k+1}
\eq\f12\sum\Sb 0<k<p^a\\k\eq r+1\,(\mo\ p-1)\endSb\bi{2k}k
-\sum\Sb 0<k<p^a\\k\eq r\,(\mo\ p-1)\endSb\bi{2k}k\ (\mo\ p^2)$$
and
$$\sum\Sb 0<k<p^a\\k\eq r\,(\mo\ p-1)\endSb C_k
\eq2\sum\Sb 0<k<p^a\\k\eq r\,(\mo\ p-1)\endSb \bi{2k}k
-\f12\sum\Sb 0<k<p^a\\k\eq r+1\,(\mo\ p-1)\endSb\bi{2k}k\ (\mo\ p^2).$$
\endproclaim
\Proof. For $k\in\N$ we have
$$\bi{2k}{k+1}+\bi{2k}k=\bi{2k+1}{k+1}=\f12\bi{2(k+1)}{k+1}.$$
Thus
$$\align&\sum\Sb 0\ls k<p^a\\k\eq r\,(\mo\ p-1)\endSb\bi{2k}{k+1}+\sum\Sb 0\ls k<p^a\\ k\eq r\,(\mo\ p-1)\endSb\bi{2k}k
\\=&\f12\sum\Sb 0\ls k<p^a\\ k\eq r\,(\mo\ p-1)\endSb\bi{2(k+1)}{k+1}
=\f12\sum\Sb 1\ls k\ls p^a\\ k\eq r+1\,(\mo\ p-1)\endSb\bi{2k}{k}
\\=&\f12\sum\Sb 0<k<p^a\\ k\eq r+1\,(\mo\ p-1)\endSb\bi{2k}{k}+R\ (\mo\ p^2)
\endalign$$
where
$$R=\cases\f12\bi{2p^a}{p^a}=\bi{2p^a-1}{p^a-1}\eq1\ (\mo\ p^2)&\t{if}\ p-1\mid r,
\\0&\t{otherwise}.\endcases$$
Therefore the first congruence in Lemma 3.2 holds. This implies
the second congruence in Lemma 3.2.
We are done. \qed

\proclaim{Lemma 3.3} Let $m,n\in\N$. Then
$$\sum_{k=0}^n\bi mk(-4)^k\bi{2(n-k)}{n-k}=4^n\prod_{k=1}^n\l(1-\f{2m+1}{2k}\r).$$
\endproclaim
\Proof. For any $k\in\N$, clearly
$$\bi{2k}k=(-4)^k\bi{-1/2}k.$$
So we have
$$\align&\sum_{k=0}^n\bi mk(-4)^k\bi{2(n-k)}{n-k}=(-4)^n\sum_{k=0}^n\bi mk\bi{-1/2}{n-k}
\\=&(-4)^n\bi{m-1/2}n=(-2)^n\prod_{k=1}^n\f{2m-2k+1}k.
\endalign$$
Therefore the desired congruence holds. \qed

\proclaim{Lemma 3.4} Let $p$ be an odd prime and let $r\in\{1,\ldots,(p-1)/2\}$.
Then
$$\align&\sum_{j=r}^{(p-1)/2}\bi{(p-1)/2}j(-4)^j\bi{2(p-1+r-j)}{p-1+r-j}
\\&\quad\eq-p\sum_{s=0}^{(p-1)/2-r}\f{\bi{2r+2s}{r+s}}{(2s+1)\bi{2s}s}\ (\mo\ p^2).
\endalign$$
\endproclaim
\Proof. If $r\ls j\ls(p-1)/2$, then $0\ls j-r<(p-1)/2$. When $s\in\N$ and $s<(p-1)/2$, clearly
$$\align&\bi{2(p-1-s)}{p-1-s}=\prod_{0<t<p-s}\f{p-1-s+t}t
\\=&\f p{s+1}\prod_{0<t\ls s}\f{p+t-s-1}t\times\prod_{s+1<t<p-s}\f{p+t-s-1}t
\\\eq&\f p{s+1}(-1)^s\prod_{0<t\ls s}\f{s-t+1}t\times\f{(p-2(s+1))!}{(p-1-s)!/(s+1)!}
\\\eq&\f{p}{s+1}(-1)^s\f{(s+1)!}{\prod_{k=s+1}^{2s+1}(p-k)}
\\\eq&\f{p(-1)^ss!}{(-1)^{s+1}\prod_{k=s+1}^{2s+1}k}
=-\f p{(2s+1)\bi{2s}s}\ (\mo\ p^2).
\endalign$$
Therefore
$$\align&\sum_{r\ls j\ls(p-1)/2}\bi{(p-1)/2}j(-4)^j\bi{2(p-1+r-j)}{p-1+r-j}
\\\eq&-p\sum_{r\ls j\ls(p-1)/2}\f{\bi{-1/2}j(-4)^j}{(2(j-r)+1)\bi{2(j-r)}{j-r}}
\\\eq&-p\sum_{r\ls j\ls(p-1)/2}\f{\bi{2j}j}{(2(r-j)+1)\bi{2(j-r)}{j-r}}\ (\mo\ p^2)
\endalign$$
and hence the desired result follows. \qed

\proclaim{Lemma 3.5} Let $p$ be an odd prime and let $a\in\Z^+$ be
even. Let $m$ be an integer not divisible by $p$ and set
$\Delta=m(m-4)$. Then
$$\sum_{k=0}^{p^a-1}\f{\bi{2k}k}{m^k}\eq
\Delta^{(p-1)/2}\sum_{k=0}^{p-1}\f{\bi{2k}k}{m^k}+\f{\da_m-\Delta^{p-1}}2\
(\mo\ p^2),$$
where $\da_m$ takes $0$ or $1$ according as  $m\eq4\ (\mo\ p)$ or
not.
\endproclaim
\Proof. By Theorem 1.1, $\sum_{k=0}^{p-1}\bi{2k}k/{m^k}\eq(\f{\Delta}p)\ (\mo\ p)$ and
$$\align\sum_{k=0}^{p^a-1}\f{\bi{2k}k}{m^k}\eq&\l(\f{\Delta}{p^{a-1}}\r)\sum_{k=0}^{p-1}\f{\bi{2k}k}{m^k}
=\l(\f{\Delta}p\r)\(\sum_{k=0}^{p-1}\f{\bi{2k}k}{m^k}-\l(\f{\Delta}p\r)\)+\l(\f{\Delta}p\r)^2
\\\eq&\Delta^{(p-1)/2}\sum_{k=0}^{p-1}\f{\bi{2k}k}{m^k}-\l(\f{\Delta}p\r)\l(\Delta^{(p-1)/2}-\l(\f{\Delta}p\r)\r)
\ (\mo\ p^2).
\endalign$$
Since
$$\align\Delta^{p-1}-\da_m=&\l(\Delta^{(p-1)/2}+\l(\f{\Delta}p\r)\r)\l(\Delta^{(p-1)/2}-\l(\f{\Delta}p\r)\r)
\\\eq&2\l(\f{\Delta}p\r)\l(\Delta^{(p-1)/2}-\l(\f{\Delta}p\r)\r)\ (\mo\ p^2),\endalign$$
the desired congruence follows from the above. \qed

\medskip
\noindent{\it Proof of Theorem 1.2}.
Let $d\in\{0,1\}$. In view of Lemma 3.1, we have
$$\align &(p-1)\sum\Sb 0<k<p^a\\k\eq r\,(\mo\ p-1)\endSb\bi{2k}{k+d}
\\\eq&\sum_{k=1}^{p^a-1}\bi{2k}{k+d}\sum_{m=1}^{p-1}m^{p(r-k)}
=\sum_{m=1}^{p-1}m^{pr}\sum_{k=1}^{p^a-1}\f{\bi{2k}{k+d}}{m^{pk}}
\ (\mo\ p^2)
\endalign$$
and hence
$$\sum\Sb 0<k<p^a\\k\eq r\,(\mo\ p-1)\endSb\bi{2k}{k+d}\ \mo\ p^2$$
only depends on the parity of $a$ by Theorem 1.1.

(i) If $a$ is odd and $r\in\{1,\ldots,p-1\}$, then by the above we have
$$\sum\Sb 0<k<p^a\\k\eq r\,(\mo\ p-1)\endSb\bi{2k}{k+d}
\eq\sum\Sb 0<k<p\\k\eq r\,(\mo\ p-1)\endSb\bi{2k}{k+d}=\bi{2r}{r+d}\ (\mo\ p^2)$$
for $d=0,1$, therefore both (1.7) and (1.8) are valid.

(ii) Now we handle the case $2\mid a$. By Lemma 3.2 it suffices to prove (1.9) for any given
$r\in\{1,\ldots,p\}$.

In light of Lemmas 3.1 and 3.5,
$$\align&(p-1)\sum\Sb 0\ls k<p^a\\k\eq r\,(\mo\ p-1)\endSb\bi{2k}k
\\\eq&\sum_{k=0}^{p^a-1}\bi {2k}k\sum_{m=1}^{p-1}m^{p(r-k)}
=\sum_{m=1}^{p-1}m^{pr}\sum_{k=0}^{p^a-1}\f{\bi{2k}k}{m^{pk}}
\\\eq&\sum_{m=1}^{p-1}m^{pr}\l(m^p(m^p-4)\r)^{(p-1)/2}\sum_{k=0}^{p-1}\f{\bi{2k}k}{m^{pk}}
\\&+\sum_{m=1}^{p-1}m^{pr}\f{\da_m-(m^p(m^p-4))^{p-1}}2\ (\mo\ p^2),
\endalign$$
where $\da_m$ is as in Lemma 3.5. (Note that $\da_{m^p}=\da_m$ since $m^p\eq m\ (\mo\ p)$.)

Observe that
$$\align&\sum_{m=1}^{p-1}m^{pr}(m^p(m^p-4))^{(p-1)/2}\sum_{k=0}^{p-1}\f{\bi{2k}{k}}{m^{pk}}
\\=&\sum_{k=0}^{p-1}\bi{2k}{k}\sum_{m=1}^{p-1}m^{p((p-1)/2+r-k)}\sum_{j=0}^{(p-1)/2}\bi{(p-1)/2}j(-4)^jm^{p((p-1)/2-j)}
\\\eq&\sum_{j=0}^{(p-1)/2}\bi{(p-1)/2}j(-4)^j\sum_{k=0}^{p-1}\bi{2k}{k}\sum_{m=1}^{p-1}m^{p(r-j-k)}\ (\mo\ p^2).
\endalign$$
So, with the help of Lemma 3.1  we have
$$\align&\f1{p-1}\sum_{m=1}^{p-1}m^{pr}(m^p(m^p-4))^{(p-1)/2}\sum_{k=0}^{p-1}\f{\bi{2k}{k}}{m^{pk}}
\\\eq&\sum_{j=0}^{(p-1)/2}\bi{(p-1)/2}j(-4)^j\sum^{p-1}\Sb k=0\\p-1\mid k+j-r\endSb\bi{2k}{k}
\\\eq&\sum_{j=0}^r\bi{(p-1)/2}j(-4)^j\bi{2(r-j)}{r-j}
\\&+\da_{r,p-1}\bi{(p-1)/2}0(-4)^0+\da_{r,p}\bi{(p-1)/2}1(-4)
\\&+\da_{r,p}\bi{(p-1)/2}0\(\bi{2\times1}{1}-\bi{2p}{p}\)
\\&+\sum_{r\ls j\ls(p-1)/2}\bi{(p-1)/2}j(-4)^j\bi{2(p-1+r-j)}{p-1+r-j}\ (\mo\ p^2).
\endalign$$
Note that
$$\bi 2{1}-\bi{2p}{p}\eq 2-2\bi{2p-1}{p-1}\eq0\ (\mo\ p^2).$$
By Lemma 3.3,
$$\align&\sum_{j=0}^r\bi{(p-1)/2}j(-4)^r\bi{2(r-j)}{r-j}=4^r\prod_{0<k\ls r}\l(1-\f p{2k}\r)
\\\eq&\cases 4^r(1-pH_r/2)\ (\mo\ p^2)&\t{if}\ 1\ls r<p-1,
\\4^r\ (\mo\ p^2)&\t{if}\ r=p-1,
\\4^r(1-pH_{p-1}/2)/2\eq 4^r/2\ (\mo\ p^2)&\t{if}\ r=p,
\endcases\endalign$$
where $H_r$  denotes the harmonic sum $\sum_{0<k\ls r}1/k$ and we note that
$$H_{p-1}=\f12\sum_{k=1}^{p-1}\l(\f1k+\f1{p-k}\r)=\f12\sum_{k=1}^{p-1}\f p{k(p-k)}\eq0\ (\mo\ p).$$
Combining the above and Lemma 3.4, we get
$$\align&\f1{p-1}\sum_{m=1}^{p-1}m^{pr}(m^p(m^p-4))^{(p-1)/2}\sum_{k=0}^{p-1}\f{\bi{2k}{k+d}}{m^{pk}}
\\\eq&-pR_p(r)+\cases4^r(1-pH_r/2)\ (\mo\ p^2)&\t{if}\ 1\ls r<p-1,
\\4^r+1\ (\mo\ p^2)&\t{if}\ r=p-1,
\\4^r/2-2p+2\ (\mo\ p^2)&\t{if}\ r=p.
\endcases
\endalign$$

Note also that
$$\align&\f1{p-1}\sum_{m=1}^{p-1}m^{pr}\l(\da_m-(m^p(m^p-4))^{p-1}\r)
\\\eq&\f1{p-1}\sum_{m=1}^{p-1}m^{pr}-\f{4^{pr}}{p-1}
-\sum_{k=0}^{p-1}\bi{p-1}k(-4)^k\f1{p-1}\sum_{m=1}^{p-1}m^{p(p-1-k+r)}
\\\eq&\da_{r,p-1}+(p+1)4^{pr}-\bi{p-1}r(-4)^r-\da_{r,p}\bi{p-1}1(-4)
-\da_{r,p-1}\bi{p-1}0
\\\eq&4^{pr}+p4^r+4(p-1)\da_{r,p}-\cases
4^r(1-pH_r)\ (\mo\ p^2)&\t{if}\ 1\ls r<p-1,
\\4^r\ (\mo\ p^2)&\t{if}\ r=p-1,\\0\ (\mo\ p^2)&\t{if}\ r=p.
\endcases
\endalign$$
So, from the above, we finally obtain
$$\sum\Sb 0\ls k<p^a\\k\eq r\,(\mo\ p-1)\endSb\bi{2k}k\eq\f{(p+1)4^r+4^{pr}}2-pR_p(r)+\da_{r,p-1}\ (\mo\ p^2).
$$
Hence
$$\sum\Sb 0<k<p^a\\k\eq r\,(\mo\ p-1)\endSb\bi{2k}k\eq\f{(p+1)4^r+4^{pr}}2-pR_p(r)\ (\mo\ p^2),$$
which is equivalent to (1.9) since
$$4^{pr}-4^r=4^r\l((1+(2^{p-1}-1))^{2r}-1\r)\eq 4^r\times 2r(2^{p-1}-1)\ (\mo\ p^2).$$

So far we have completed the proof of Theorem 1.2. \qed

\medskip
\noindent{\it Proof of Corollary 1.2}.
Recall that $H_{p-1}\eq0\ (\mo\ p)$. As observed by Eisenstein,
$$\align&\f{2^p-2}p=\sum_{k=1}^{p-1}\f1p\bi{p}k=\sum_{k=1}^{p-1}\f1k\bi{p-1}{k-1}
\\\eq&\sum_{k=1}^{p-1}\f{(-1)^{k-1}}k\eq\sum_{k=1}^{p-1}\f{(-1)^{k-1}-1}k
=\sum_{j=1}^{(p-1)/2}\f{-2}{2j}=-H_{(p-1)/2}\ (\mo\ p).
\endalign$$
It is easy to see that
$$\align C_{p-1}=&\f1{p-1}\bi{2p-2}{p-2}=\f1{2p-1}\prod_{k=1}^{p-1}\l(1+\f pk\r)
\\\eq&-(1+2p)(1+pH_{p-1})\eq-1-2p\ (\mo\ p^2)
\endalign$$
and
$$\align C_{(p-1)/2}=&\f2{p+1}\bi{p-1}{(p-1)/2}
\\=&\f2{p+1}(-1)^{(p-1)/2}\prod_{k=1}^{(p-1)/2}\l(1-\f pk\r)
\\\eq&2(1-p)(-1)^{(p-1)/2}(1-pH_{(p-1)/2})
\\\eq&2(-1)^{(p-1)/2}(1-p-pH_{(p-1)/2})
\\\eq& 2(-1)^{(p-1)/2}(2^p-p-1)\ (\mo\ p^2).
\endalign$$
So, by Theorem 1.2(i), (1.13)-(1.15) hold in the case $2\nmid a$.

From now on we assume that $a$ is even.

Applying (1.12) with $r=p-1$ we immediately get (1.13).
As
$$R_p\l(\f{p-1}2\r)=\bi{p-1}{(p-1)/2}\eq(-1)^{(p-1)/2}\ (\mo\ p)$$
and  $R_p((p+1)/2)=0$, by (1.11) we have
$$\align\sum\Sb 0<k<p^a\\k\eq(p-1)/2\,(\mo\ p-1)\endSb C_k
\eq& 4^{(p-1)/2}(2-2^p)-p2(-1)^{(p-1)/2}
\\\eq&2-2^p+(-1)^{(p+1)/2}2p\ \ (\mo\ p^2).
\endalign$$
This proves (1.15).

To obtain (1.14) we need to compute $R_p(1)$ and $R_p(2)$ modulo $p$.
Observe that
$$\align R_p(1)=&\sum_{s=0}^{(p-1)/2-1}\f{2\bi{2s+1}{s}}{(2s+1)\bi{2s}s}=\sum_{s=0}^{(p-3)/2}\f2{s+1}
\\=&2H_{(p-1)/2}\eq2\times\f{2-2^p}p\ (\mo\ p).
\endalign$$
When $p\gs 5$, we have
$$\align R_p(2)=&\sum_{s=0}^{(p-1)/2-2}\f{2\bi{2s+3}{s+1}}{(2s+1)\bi{2s}s}=\sum_{s=0}^{(p-5)/2}\f{4(2s+3)}{(s+1)(s+2)}
\\=&4\sum_{s=0}^{(p-5)/2}\l(\f1{s+1}+\f1{s+2}\r)=4(H_{(p-3)/2}+H_{(p-1)/2}-1)
\\=&8H_{(p-1)/2}-4\l(\f2{p-1}+1\r)
\eq 8\times\f{2-2^p}p+4\ (\mo\ p).
\endalign$$
In the case $p=3$, as $R_3(2)=0$ we also have $R_p(2)\eq 8(2-2^p)/p+4\ (\mo\ p)$.
Applying (1.11) with $r=1$, we obtain
$$\align\sum\Sb 0<k<p^a\\k\eq1\,(\mo\ p-1)\endSb C_k&\eq 4(2-2^p)-p\l(2R_p(1)-\f{R_p(2)}2\r)
\\\eq& 4(2-2^p)-p(-2)\ (\mo\ p^2).
\endalign$$
So (1.14) follows.

The proof of Corollary 1.2 is now complete. \qed


\bigskip

 \widestnumber\key{ZPS}

 \Refs

\ref\key CP\by R. Crandall and C. Pomerance
\book Prime Numbers: A Computational Perspective
\publ 2nd edition, Springer, New York, 2005\endref

\ref\key GKP\by R. L. Graham, D. E. Knuth and O. Patashnik
 \book Concrete Mathematics\publ 2nd ed., Addison-Wesley, New York\yr 1994\endref

\ref\key J\by W. J\"anichen\paper \"Uber die Verallgemeinerung einer Gauss'schen Formel aus
der Theorie der h\"oheren Kongruenzen\jour Sitzungsber. Berl. Math. Ges.\vol 20\yr 1921\pages 23--29\endref

\ref\key PS\by H. Pan and Z. W. Sun\paper A combinatorial identity
with application to Catalan numbers \jour Discrete Math.\vol
306\yr 2006\pages 1921--1940\endref

\ref\key Sm\by C. J. Smyth\paper A coloring proof of a generalization of Fermat's little theorem
\jour Amer. Math. Monthly\vol 93\yr 1986\pages 469--471\endref

\ref\key St\by R. P. Stanley\book Enumerative Combinatorics \publ
Vol. 2, Cambridge Univ. Press, Cambridge, 1999\endref

\ref\key SS\by Z. H. Sun and Z. W. Sun\paper Fibonacci numbers and Fermat's last theorem
\jour Acta Arith.\vol 60\yr 1992\pages 371--388\endref

\ref\key S92\by Z. W. Sun\paper Reduction of unknowns in Diophantine representations
\jour Sci. China Ser. A \vol 35\yr 1992\pages no.3, 257--269\endref

\ref\key S02\by Z. W. Sun\paper On the sum $\sum_{k\eq r\,(\mo\ m)}\bi nk$
and related congruences\jour Israel J. Math.
\vol 128\yr 2002\pages 135--156\endref

\ref\key S06\by Z. W. Sun\paper Binomial coefficients and quadratic fields
\jour Proc. Amer. Math. Soc.\vol 134\yr 2006\pages 2213--2222\endref

\ref\key S09\by Z. W. Sun\paper Various congruences involving binomial coefficients and higher-order Catalan numbers
\jour preprint, arXiv:0909.3808. {\tt http://arxiv.org/abs/0909.3808}\endref


\ref\key ST1\by Z. W. Sun and R. Tauraso\paper New congruences for
central binomial coefficients \jour Adv. in Math. \vol 45\yr
2010\pages 125--148\endref

\ref\key ST2\by Z. W. Sun and R. Tauraso\paper On some new
congruences for binomial coefficients \jour Int. J. Number Theory\pages in press.
{\tt http://arxiv.org/abs/0709.1665}\endref


\ref\key V\by E. B. Vinberg\paper On some number-theoretic conjectures of V. Arnold
\jour Jpn. J. Math.\vol 2\yr 2007\pages 297--302
\endref

\ref\key ZPS\by L. L. Zhao, H. Pan and Z. W. Sun\paper Some congruences for the second-order Catalan numbers
\jour Proc. Amer. Math. Soc.\vol 138\yr 2010\pages 37--46\endref

\endRefs
\enddocument